\documentclass{amsart}
\usepackage{amssymb,latexsym}
\usepackage[all]{xy}
\theoremstyle{plain}
\newtheorem*{main*}{Main Theorem}
\newtheorem{theorem}{Theorem}

\newtheorem{lemma}[theorem]{Lemma}
\theoremstyle{definition}

\newtheorem{remark}[theorem]{Remark}
\newtheorem{question}[theorem]{Question}
\newtheorem{claim}[theorem]{Claim}

\DeclareMathOperator{\hh}{h}

\DeclareMathOperator{\Res}{Res}
\DeclareMathOperator{\len}{length}
\newcommand{\I}{\mathcal{I}}
\newcommand{\A}{\sf{{A}}}
\newcommand{\B}{\sf{B}}
\newcommand{\C}{\sf{C}}
\newcommand{\PP}{\mathbf P}
\usepackage[naturalnames=true]{hyperref}

\date{}

\begin{document}

\title[Quartuple points]{Postulation of general quartuple \\
fat point schemes in $\mathbf{P}^3$}

\author{Edoardo Ballico}
\address{Dept. of Mathematics\\
 University of Trento\\
38050 Povo (TN), Italy}
\email{ballico@science.unitn.it}

\author{Maria Chiara Brambilla}
\address{Dept. of Mathematics ``G.~Castelnuovo'', University of Rome
  ``La Sapienza", Piazzale Aldo Moro 2, 00185, Rome, Italy}
\email{brambilla@math.unifi.it, brambilla@mat.uniroma1.it}

\thanks{Both authors were partially supported by MIUR and GNSAGA of
  INdAM (Italy).}

\maketitle

\begin{abstract}
We study the postulation of a general union $Y$ of double, triple, and quartuple points of $\PP^3$.
We prove that $Y$ has the expected postulation in degree $d\ge 41$, using the Horace differential
lemma. We also discuss the cases of low degree with the aid of computer algebra.
\end{abstract}

\section{Introduction}

In this paper we study the postulation of general fat point
schemes of $\PP^3$. A {\em fat point} $mP$ is a zero-dimensional
subscheme of $\PP^3$ supported at a point $P$ and with $({\mathcal
{I}}_{P,\PP^3})^m$ as its ideal sheaf. A {\em general fat point
scheme} $Y=m_1P_1+\ldots+m_kP_k$, with $m_1\ge\ldots\ge m_k\ge1$,
is a general zero-dimensional scheme such that $Y_{\textrm{red}}$
is a union of $k$ points and for each $i$ the connected component
of $Y$ supported at $P_i$ is the fat point $m_iP_i$. We call
multiplicity of $Y$ the maximal multiplicity, $m_1$, of its
components. We recall that $\len(mP)={{m+2}\choose{3}}$, for any
$m\ge1$.

Studying the postulation of $Y$ means to compute the dimension of the space
of hypersurfaces of any degree containing the scheme $Y$. In other
words this is equivalent to
compute the dimension of the space of homogeneous polynomials of any degree
vanishing at the point $P_i$ and with all their derivatives,
up to multiplicity $m_i-1$, vanishing at $P_i$.
We say that $Y$ has {\em good postulation} if such a dimension is the expected one.

This problem has been investigated by many authors in the case of $\PP^2$.
In particular we recall the important Harbourne-Hirschowitz conjecture
(see the survey \cite{ciliberto} and references therein). This conjecture characterizes all the
general fat point schemes not having good postulation, and has been
proved in some special cases.
We mention also an analogous conjecture in the case of $\PP^3$,
due to Laface and Ugaglia (see \cite{laface-ugaglia}).
In the case of general unions of double points, that is when $m_i=2$ for any $i$,
the famous Alexander-Hirschowitz theorem gives a
complete answer in the case of $\PP^n$, for any $n\ge 2$,
(see  \cite{AHinv,AHjag}, for a survey see \cite{BO}).
For arbitrary multiplicities and arbitrary projective
varieties there is a beautiful asymptotic theorem by
Alexander and Hirschowitz \cite{AH}.

Here we will study the case of general fat point scheme
$Y\subset\PP^3$ of multiplicity $4$. The case of multiplicity $3$
was considered by the first author in \cite{ballico-tripli}, where
he proved that a general union $Y \subset \PP^3$ of triple and
double points has good postulation in degree $d\ge7$.
\newpage
Our main result is the following:

\begin{theorem}\label{i1}
Assume $\textrm{char}(\mathbf{K})\neq 2,3$.
Fix non-negative integers $d, x, y, z$ such that $d \ge 41$.
Let $Y \subset {\bf {P}}^3$
be a general union of $x$ 4-points, $y$ 3-points and $z$
2-points. Then $Y$ has good postulation, i.e.
\begin{itemize}
\item
if $20x + 10y+4z \le \binom{d+3}{3}$,
then $\hh^1({\bf{P}}^3,\mathcal {I}_Y(d)) = 0$,
\item if $20x+10y + 4z \ge \binom{d+3}{3}$, then
$\hh^0({\bf {P}}^3,\mathcal {I}_Y(d))= 0$.
\end{itemize}
\end{theorem}

The proof is based on the well known {\em Horace differential lemma}.
We point out that this asymptotic result is not proved by induction on the degree,
hence it does not depend on the cases of low degree.

The cases where $d\le 40$ can be analyzed with the help of
computer algebra. We have checked that if $d\le 8$ there exist some
cases where a general fat point scheme $Y$ of multiplicity $4$
does not have good postulation in degree $d$.
This happens in particular if
the number of quartuple points contained in $Y$ is high.
On the other hand, we found that if
$9\le d \le 13$ any general fat point scheme $Y$ of multiplicity $4$
has good postulation in degree $d$.
We expect that the same is true also for $14\le d\le 40$, even if
we did not perform the computations.

With the same kind of computation one may start to investigate the cases of fat
point schemes of multiplicity higher than $4$ for low degree. In
Section 4.1 we have collected some partial results in this direction.

These numerical experiments lead us to pose the following question,
which we believe is interesting even for low multiplicities cases:
\begin{question}
Let $Y \subset \PP^n$ be a general
fat point scheme of multiplicity $m\ge2$.
Let $d(n,m)$ be a function such that for any $d\ge d(n,m)$ the scheme
$Y$ has good postulation in degree $d$.
For fixed $n$ is it possible to take as $d(n,m)$ a function polynomial
(or even linear) in $m$?
Is it possible to take $d(3,m)=3m$?
\end{question}

Note that by \cite[Example 7.7]{laface-ugaglia} we know that
$d(3,m)>2m$. We also know that $d(3,m)>2m+1$.
In fact, the referee suggested us the following example:
$9$ general $9$-points of $\PP^3$ have not good postulation in degree $19$.

Notice that our question concerns an upper estimate which is not sharp.
It seems difficult to find a sharp estimate, and of course it would be very interesting.
For other results related to this subject see also \cite{volder-laface}
and \cite{dumnicki}.

\medskip

Here is the plan of the paper.
In Section 2 we give some preliminary lemmas.
Section 3 is devoted to the proof of the main result of the paper
(Theorem \ref{i1}), while in Section 4 we give some details on the
cases of low degree.

\section{Preliminaries}
Throughout the paper we will work on the $n$-dimensional projective
space $\PP^n$ over an algebraically closed field $\mathbf{K}$.

For any smooth $n$-dimensional connected projective variety $A$,
any $P\in A$ and any integer $m>0$, an {\em $m$-fat point of $A$} (or just
$m$-point) $\{mP,A\}$ is defined to be the $(m-1)$-th
infinitesimal neighborhood of $P$ in $A$, i.e.\ the closed subscheme
of $A$ with $(\mathcal {I}_{P,A})^m$ as its ideal sheaf.
Thus $\{mP,A\}_{red} = \{P\}$ and $\len(\{mP,A\})= \binom{n+m-1}{n}$.
We will write $mP$ instead of $\{mP,A\}$ when the space $A$ is clear
from the context, and mostly we will consider $A={\bf {P}}^n$ for $n=2,3$.
We call {\em general fat point scheme of $A$} a union $Y=m_1P_1+\ldots+m_kP_k$, with
$m_1\ge\ldots\ge m_k\ge1$, and $P1,\ldots,P_k$ general points of
$\PP^n$. We denote $\deg(Y)=\sum \len(m_iP_i)$.

Given a positive integer $d$, we will say that a zero-dimensional
scheme $Y$ of $\PP^n$ has {\em good postulation in degree $d$} if the following
conditions hold:
\begin{enumerate}
\item[(a)]
if $\deg(Y) \le \binom{n+d}{n}$, then
$\hh^1(\PP^n,\mathcal{I}_{Y}(d))=0,$
\item[(b)] if $\deg(Y) \ge\binom{n+d}{n}$, then
$\hh^0(\PP^n,\mathcal {I}_{Y}(d))=0.$
\end{enumerate}

Given a general fat point scheme $Y$ of $\PP^n$ and a hyperplane
$H\subset\PP^n$ we will call {\em trace} of $Y$ the subscheme $Y\cap
H\subset H$ and {\em residual} of $Y$ the scheme $\Res_H(Y)\subset\PP^n$ with
ideal sheaf $\I_Y:\mathcal{O}_{\PP^n}(-H)$.
Notice that if $X$ is a $m$-point supported on $H$, then its trace
$X\cap H$ is a $m$-point of $H$ and its residual $\Res_H(X)$ is a
$(m-1)$-point of $\PP^n$.

The trace and the residual of a fat point scheme $Y$ of $\PP^n$ fit in the
following well known {\em Castelnuovo exact sequence}
$$0\to \I_{\Res_H(Y)}(d-1)\to \I_Y(d)\to \I_{Y\cap H}(d)\to 0.$$
A straightforward consequence of the Castelnuovo exact sequence is the
following form of the so called {\em Horace lemma},
which we will often use in the sequel. For more details see e.g. \cite[Section 4]{BO}.

\begin{lemma}\label{castelnuovo}
Let $H\subset \PP^n$ be a hyperplane and $Y\subset \PP^n$ a fat point scheme of $\PP^n$.
Then we have
$$ \hh^0(\PP^n,\I_Y(d))\le
\hh^0(\PP^n,\I_{\Res_H(Y)}(d-1))+\hh^0(H,\I_{Y\cap H}(d))$$
$$ \hh^1(\PP^n,\I_Y(d))\le
\hh^1(\PP^n,\I_{\Res_H(Y)}(d-1))+\hh^1(H,\I_{Y\cap H}(d))$$
\end{lemma}

The basic tool we will need is the so called {\em Horace differential
lemma}, introduced by Alexander and Hirschowitz. This technique allows us to take a {\em differential trace}
and a {\em differential residual}, instead of the classical ones.
For an explanation of the geometric intuition of the Horace
differential lemma see \cite[Section 2.1]{AH}. Here we give only an
idea of how the lemma works.

Let $Y$ be an $m$-point of $\PP^n$ supported on a hyperplane
$H\subset\PP^n$. Following the language of Alexander and Hirschowitz
we can describe $Y$ as formed by infinitesimally piling up some
subschemes of $H$, called {\em layers}.
For example the layers of a $3$-point $\{3P,\PP^n\}$ are $\{3P,H\}$,
$\{2P,H\}$, and $\{P,H\}$.
Then the differential trace can be any of these layers and the
differential residual is a {\em virtual} zero-dimensional scheme
formed by the remaining layers.
We will denote these virtual schemes by writing the subsequent layers
from which they are formed.
These layers are obtained intersecting with the hyperplane $H$
and taking the residual many times.
In particular the notation e.g.\ $X= (\{3P,H\},  \{2P,H\})$
means that $X\cap H= \{3P,H\}$ and $\Res_H(X)\cap H= \{2P,H\}$, and,
finally, $\Res_H(\Res_H(X))\cap H= \emptyset$.

In this paper we will apply several times the following result which
is a particular case of the Horace differential lemma (see \cite[Lemma
2.3]{AH}).
\begin{lemma}[Alexander-Hirschowitz]
\label{alehir}
Fix an integer $m\ge 2$ and
assume that $\textrm{char}(\mathbf{K})=0$
or $\textrm{char}(\mathbf{K})>m$.
Let $X$ be an $m$-point of $\PP^n$ supported at $P$
and $H\subset\PP^n$ a hyperplane.
Then for $i=0,1$ we have
$$ \hh^i(\PP^n,\I_X(d))\le \hh^i(\PP^n,\I_{R}(d-1))+\hh^i(H,\I_{T}(d))$$
where the {\em differential residual } $R$ and the {\em differential trace} $T$
are virtual scheme of the following type:
\begin{center}
\begin{tabular}{lllll}
$(i)$ &  $m=2$: & $T=\{P,H\}$;& $R=\{2P,H\}$ &$(1,3)$\\
$(ii)$& $m=3$: & $T=\{P,H\}$;& $R= (\{3P,H\},  \{2P,H\})$&$(1,6,3)$\\
$(iii)$& $m=3$:& $T=\{2P,H\}$;& $R= (\{3P,H\},  \{P,H\})$&$(3,6,1)$\\
$(iv)$&$m=4$:& $T=\{P,H\}$ ;& $R=( \{4P,H\}, \{3P,H\}, \{2P,H\})$&$(1,10,6,3)$\\
$(v)$&$m=4$: & $T=\{2P,H\}$; & $R=(\{4P,H\}, \{3P,H\}, \{P,H\})$&$(3,10,6,1)$\\
$(vi)$&$m=4$: & $T=\{3P,H\}$; &$R=( \{4P,H\},  \{2P,H\}, \{P,H\})$&$(6,10,3,1)$
\end{tabular}
\end{center}
\end{lemma}

In the previous lemma, for each case in the statement we write in the last column
the list of the lengths of the fat points of $H$ that we will obtain
intersecting many times with $H$.
Throughout the paper, when we will apply Lemma \ref{alehir}, we will
specify which case we are considering by recalling this sequence of the
lengths.
For example if we apply Lemma \ref{alehir}, case $(i)$, we
will say the we apply the lemma with respect to the sequence $(1,3)$.

\begin{remark}\label{facile}
Let $X \subseteq Y \subset {\bf {P}}^n$ zero-dimensional schemes.
Then it is immediate to see that $\hh^0({\bf {P}}^n,\mathcal {I}_Y(d)) \le \hh^0({\bf
  {P}}^n,\mathcal {I}_X(d))$.
\end{remark}

We recall here a particular case of a result of Mignon
(see \cite[Theorem 1]{m}).
\begin{lemma}[Mignon]
\label{mignon} Let $X\subset \PP^2$ be a general fat point scheme of
multiplicity $4$ (that is a general collection of multiple points of
multiplicity at most $4$) and $d\ge 12$. Then $X$ has good postulation, i.e. we have
\begin{enumerate}
\item[(a)]
if $\deg(X) \le \binom{d+2}{2}$, then
$\hh^1(\PP^2,\mathcal{I}_{X}(d))=0,$
\item[(b)] if $\deg(X) \ge\binom{d+2}{2}$, then
$\hh^0(\PP^2,\mathcal {I}_{X}(d))=0.$
\end{enumerate}
If $X\subset \PP^2$ is a general fat point scheme of
multiplicity $3$, and $d\ge9$, then $X$ has good postulation.
\end{lemma}

The following lemma is equivalent to \cite[Remark 2]{ballico-tripli}.
We give here a complete proof for the reader's convenience.
\begin{lemma}\label{remark2}
Fix integers $d>0$, $z>0$, $\gamma\ge0$,  a hyperplane $H\subset {\bf {P}}^n$
and a zero-dimensional scheme $Y\subset {\bf {P}}^n$.
Let $X$ be the union of $Y$ and $z$ general simple points supported on $H$.
If the following conditions
\begin{equation}\label{G}
\hh^0({\bf {P}}^n,\I_Y(d))\le\gamma+z, \quad\mbox{and}
\quad \hh^0({\bf {P}}^n,\I_{\Res_H(Y) }(d-1))\le\gamma,
\end{equation}
take place, then it follows that
$$
\hh^0({\bf {P}}^n,\I_X(d))\le\gamma.
$$
Equivalently if the following conditions
$$\hh^1(\PP^n,\I_Y(d))\le \max(0,\gamma+\deg(X)-\binom{d+n}{n})=:\beta,$$
and
$$\hh^1(\PP^n,\I_{\Res_H(Y)}(d-1))\le
\max(0,\gamma+\deg(\Res_H(Y))-\binom{d+n-1}{n}),$$
take place, then it follows that
$$
\hh^1(\PP^n,\I_X(d))\le \beta.
$$
\end{lemma}

\begin{proof}
Notice that, since for any scheme $Z$ and any integer
$d$ we have
$$\hh^1({\bf {P}}^n,\I_Z(d))=\hh^0({\bf
  {P}}^n,\I_Z(d))-\binom{d+n}{n}+\deg(Z),$$
then the two formulations of the lemma are equivalent, since clearly
we have $\deg(X)=\deg(Y)+z$.

Let us assume that the two conditions in \eqref{G} hold and,
for any positive integer $p$, let us denote by $Y_p$ the union of $Y$ and
$p$ general simple points of $H$.
Let $r$ be the maximal integer $p$ such that
$\hh^0({\bf {P}}^n,\I_{Y_p}(d))=\hh^0({\bf {P}}^n,\I_{Y}(d))-p$.
Obviously $0\leq r \leq \hh^0({\bf {P}}^n,\I_{Y}(d))\leq \gamma + z
\le z$.
Since $Y_{r+1}\setminus Y_r$ is a general point of $H$, it follows that $H$
is contained in the base locus of the linear system $|\I_{Y_r}(d)|$.
This implies that
$\hh^0({\bf {P}}^n,\I_{Y_r}(d))=\hh^0({\bf
  {P}}^n,\I_{\Res_H(Y)}(d-1))\le \gamma$.
Then, since $Y_r$ can be identified with a subscheme of $X$, by Remark
\ref{facile}, we conclude that $\hh^0({\bf {P}}^n,\I_X(d))\leq \gamma$.
\end{proof}

The following numerical lemma will be used in the sequel.
\begin{lemma}\label{c1}
Fix non-negative integers $t, a,b,c,e,f,g$ such that $t \ge 14$,
\begin{equation}\label{eqc1}
10a +6b+3c +u +6e+3f +g \le \binom{t+2}{2}
\end{equation}
and $(e,f,g)$ is one of the following triples: $(0,0,0)$,
$(0,0,1)$, $(0,0,2)$, $(0,1,0)$, $(0,1,1)$, $(0,1,2)$, $(1,0,0)$,
$(1,0,1)$, $(1,0,2)$, $(1,1,0)$. Then we get the following
inequality
\begin{align}
\label{eqc2} &6a+3b+c+ 10e +10f + 10g \le \binom{t+1}{2}.
\end{align}
Moreover, if $e+f+g\le2$, then
\eqref{eqc2} holds for any $t\geq12$.
If $e=f=g=0$, then
\eqref{eqc2} holds for any $t\ge 3$.
\end{lemma}

\begin{proof}
In order to prove (\ref{eqc2}),
it is sufficient to check the inequality
$$4a +3b+2c +u -4e-7f-9g \ge t+1.$$
From (\ref{eqc1}) it follows that
$$10a \ge \binom{t+2}{2}-6b-3c -u -9,$$
hence the inequality above comes from
$$\frac{2}{5}\left(\binom{t+2}{2} -9\right) -25 \ge t+1,$$
which is true for all $t \ge 14$.
If $e+f+g\le2$, then we have to check
$$\frac{2}{5}\left(\binom{t+2}{2} -9\right) -18 \ge t+1,$$
which holds for any $t\geq12$.
Finally, the last statement follows easily.
\end{proof}

\section{Proof of the main theorem}

This section is devoted to the proof of Theorem \ref{i1}.
Throughout the section we assume that the characteristic of the base
field $\mathbf{K}$ is different from $2$ and $3$, and we fix an
hyperplane $H\subset\PP^3$.

In the different steps of the proof we will work with zero-dimensional
schemes a little more general than a union of fat points. In
particular, we will say that a zero-dimensional
scheme $Y$ is {\em of type $(\star)$} if its irreducible components
are of the following type:
\begin{itemize}
\item[-]
$m$-points, with $2\le m\le4$ supported at general points of $\PP^3$,
\item[-]
$m$-points, with $1\le m\le 4$, or virtual schemes arising as residual
in the list of Lemma
\ref{alehir}, supported at general points of $H$.
\end{itemize}

In the following lemma we describe a basic step that we will apply
several times in the sequel.

\begin{lemma}\label{lemma-a}
Let $Y$ be a zero-dimensional scheme of type $(\star)$.
For $2\le i \le 4$, let $c_i$ be the number of $i$-points of $Y$ not
supported in $H$.
If the following condition holds
\begin{equation}\label{condizione}
\beta:=\binom{t+2}{2}-\deg(Y\cap H)\geq0,
\end{equation}
then it is possible to degenerate $Y$ to a scheme $X$ such that
one of the following possibilities is verified:
\begin{enumerate}
\item[(I)] $\deg(X \cap H)=\binom{t+2}{2}$,
\item[(II)] $\deg(X \cap H)<\binom{t+2}{2}$, and all the irreducible
  components of $X$ are supported on $H$. This is possible only if
  $c_2+c_3+c_4\le2$ and $c_2+c_3+c_4<\beta$.
\end{enumerate}
In both cases we also have
\begin{equation}\label{caldo}
\deg(\Res_H(X)\cap H)\le \binom{t+1}{2}
\end{equation}
\end{lemma}

\begin{proof}
First of all we can assume that $\beta\ge0$ is minimal.
Indeed we can change the scheme $Y$ by specializing on $H$ some other
component which are not supported on $H$. Let us denote now  by $Y'$
the union of the connected components of $Y$ intersecting $H$.

By minimality of $\beta$ it follows that if $c_2>0$ then $\beta< 3$,
if $c_2=0$ and $c_3>0$ then $\beta< 6$, if $c_2=c_3=0$ and $c_4>0$ then
$\beta< 10$. If  $c_2=c_3=c_4=0$ and $\beta>0$, we are obviously in case
(II).

We degenerate now $Y$ to a scheme $X$ described as follows.
The scheme $X$ contains all the connected components of $Y'$.
Write
$$\beta= 6e+3f+g$$
 for a unique triple of non-negative integers $(e,f,g)$ in the following list:
$(0,0,0)$, $(0,0,1)$, $(0,0,2)$, $(0,1,0)$, $(0,1,1)$, $(0,1,2)$, $(1,0,0)$,
$(1,0,1)$, $(1,0,2)$,  $(1,1,0)$ (i.e. in the list of Lemma \ref{c1}).
Notice that if $c_2>0$ then $e=f=0$ and $g\le 2$, if $c_2=0$ and $c_3>0$ then
$e=0$ and $f+g\le 3$, if $c_2=c_3=0$ and $c_4>0$ then $e+f+g\le3$.

Consider first the case $c_2>0$ and recall that in this case $e=f=0$
and $g\le2$. Assume now $c_2 \ge g$. Take as $X$ a general union
of $Y'$, $c_4$ 4-points, $c_3$ 3-points, $(c_2-g)$ 2-points,
$g$ virtual schemes obtained applying Lemma \ref{alehir}
at $g$ general points of $H$ with respect to the sequence
$(1,3)$. Clearly we have $\deg(X\cap H)=\binom{t+2}{2}$ and we are in case (I).

Let us see now how to specialize $Y$ to $X$ in the remaining cases.
If $c_2 =1<g$ and $c_3+c_4\ge1$, then in the previous step we apply
Lemma \ref{alehir} using the unique 2-point and one 3-point (or
4-point respectively) with respect to the sequence
$(1,6,3)$ (or $(1,10,6,3)$ respectively) and we conclude in the
same way, obtaining (I).
If $c_2 =1<g$ and $c_3=c_4=0$, then we apply  Lemma \ref{alehir} to
the unique double point with respect to the sequence $(1,3)$, and we
are in case (II).

Assume now $c_2=0$ and $c_3>0$.
If $c_3\ge f+g$ we take as $X$ a general union
of $Y'$, $c_4 $ 4-points, $c_3-f-g$ $3$-points,
$f$ virtual schemes obtained applying Lemma \ref{alehir} at $f$
general points of $H$ with respect to
the sequence $(3,6,1)$ and $g$ virtual schemes obtained applying
Lemma \ref{alehir} at $g$ general points of $H$ with respect to the sequence
$(1,6,3)$.
If $0<c_3<f+g$ and $c_4\ge f+g-c_3$, then in the previous step we
apply Lemma \ref{alehir} using $c_3$ 3-point, and $(f+g-c_3)$ 4-points,
with respect to the sequences $(3,10,6,1)$ or $(1,10,6,3)$.
In all these cases we clearly have $\deg(X\cap H)=\binom{t+2}{2}$, so we are in case (I).
If  $c_2=0$, $0<c_3<f+g$ and $c_4< f+g-c_3$, then we have either
$c_3\le1$ and $c_4\le1$, or $c_3=2$ and $c_4=0$, and in both cases
$\beta>c_3+c_4$. In this cases we can specialize all the components on
$H$, possibly applying Lemma \ref{alehir} and we are in case (II).

Now, assume that $c_2=c_3=0$ and $c_4>0$.
If $c_4\ge e+f+g$, then
we take as $X$ a general union
of $Y'$, $(c_4-e-f-g)$ $4$-points, $e$ virtual schemes obtained applying
Lemma \ref{alehir} at $e$ general points of $H$ with respect to the
sequence $(6,10,3,1)$,
$f$ virtual schemes obtained applying Lemma \ref{alehir} at $f$
general points of $H$ with respect to the sequence $(3,10,6,1)$ and
$g$ virtual schemes obtained applying Lemma \ref{alehir} at $g$
general points of $H$ with respect to the sequence $(1,10,6,3)$.
Thus we have again $\deg(X\cap H)=\binom{t+2}{2}$, that is we are in case (I).
If $c_2=c_3=0$ and  $0<c_4<e+f+g$, then we are in case
(II), because we can specialize all the quartuple points on $H$
(possibly applying Lemma \ref{alehir}), since
$c_4\le e+f+g+1\le 2$ and $\beta>c_4$.

Finally, we note that the property \eqref{caldo} follows immediately by
the construction above and by Lemma \ref{c1}.
\end{proof}

Given a scheme $Y$ of type $(\star)$ satisfying
\eqref{condizione}, we will say that $Y$ is of type
(I) if, when we apply Lemma \ref{lemma-a} to $Y$, we are in case
(I). Otherwise we say that $Y$ is of type (II).

We fix now (and we will use throughout this section)
the following notation, for any integer $t$:
given a scheme $Y_t$ of type $(\star)$ and satisfying
\eqref{condizione}, we will denote by $X_t$ the
specialization described in Lemma \ref{lemma-a}.
We write the residual $\Res_H(X_t)=Y_{t-1}\cup Z_{t-1}$,
where $Y_{t-1}$ is the union of all unreduced components of
$\Res_H(X_t)$ and $Z_{t-1}=\Res_H(X_t)\setminus Y_{t-1}$.
Clearly  $Z_{t-1}$ is the union of finitely many simple points of
$H$.
Thus at each step $t \mapsto t-1$, we will have
$$Y_t \mapsto X_t \mapsto  \Res_H(X_t)= Y_{t-1}\cup Z_{t-1}.$$
For any integer $t$, we set $z_{t}:= \sharp (Z_{t})$,
$\alpha_{t}:=\deg(Y_{t})=\deg(X_{t})$, and
$$\delta_t:=\max\left(0,\binom{t+2}{3} -\deg(Y_{t-1}\cup
Z_{t-1})\right).$$

We fix the following statements:
\begin{itemize}
\item[-] ${\A}(t)=\{Y_t \mbox{ has good postulation in degree $t$}\}$,
\item[-] ${\B}(t)=\{\Res_H(X_t) \mbox{ has good postulation in degree $t-1$}\}$,
\item[-] ${\C}(t)=\{\hh^0({\bf {P}}^3,\mathcal {I}_{\Res_H(Y_{t-1})}(t-2)) \le \delta_t\}$.
\end{itemize}

\begin{claim}\label{claim-a0}
Fix $t\ge13$.
If $Y_t$ is a zero-dimensional scheme of type (II), then
it has good postulation, i.e. ${\A}(t)$ is true.  Moreover also ${\B}(t)$ is true.
\end{claim}

\begin{proof}
Since $Y_t$ is of type (II),
when we apply Lemma \ref{lemma-a} to $Y_t$, we obtain a
specialization $X_t$ whose all irreducible components are supported on
$H$ and such that $\deg(X_t\cap H)\leq\binom{t+1}{2}$.

We prove now the vanishing $\hh^1({\bf{P}}^3,\mathcal {I}_Y(t)) = 0$.
By semicontinuity, it is enough to prove the
vanishing
$\hh^1({\bf{P}}^3,\mathcal {I}_X(t)) = 0$.
Notice that by taking the residual of $X_t$ with respect to $H$ for at most
four times we get at the end the empty set.

Since $\deg(X_t\cap H)\le\binom{t+2}{2}$, and $t\ge12$,
by Lemma \ref{mignon} it follows the vanishing
$\hh^1({\bf{P}}^2,\mathcal{I}_{X\cap H}(t)) = 0.$
Let us denote by $R_{t-1}$ the residual $\Res_H(X)$ and recall that
any component of $R_{t-1}$ is supported on $H$.
We need to check now that
$\hh^1({\bf{P}}^3,\mathcal{I}_{R_{t-1}}(t)) = 0.$

In order to do this
we take again the trace and the residual with respect to $H$.
By \eqref{caldo} we know that
$\deg(\Res_H(X_t)\cap H)\leq\binom{t+1}{2}$
then again by Lemma \ref{mignon}, since $t-1\ge12$, we have
$\hh^1({\bf{P}}^2,\mathcal{I}_{R_{t-1}\cap H}(t-1)) = 0.$

We repeat now this step taking $R_{t-2}:=\Res_H(R_{t-1})$ and noting
that the trace $R_{t-2}\cap H$ has degree less or equal than
$\binom{t}{2}$, by Lemma \ref{c1}.
Moreover this time the scheme $R_{t-2}\cap H$ cannot contain quartuple
points,
in fact it is a general union of triple, double and simple points.
Hence by Lemma \ref{mignon}, since $t-2\ge9$ we have
$\hh^1({\bf{P}}^2,\mathcal{I}_{R_{t-2}\cap H}(t-2)) = 0.$

We repeat once again the same step and we obtain
$R_{t-3}:=\Res_H(R_{t-2})$. Now the trace ${R_{t-3}\cap H}$ contains
only double or simple points and so we have again the vanishing
$\hh^1({\bf{P}}^2,\mathcal{I}_{R_{t-3}\cap H}(t-3)) = 0,$
by the Alexander-Hirschowitz Theorem, since $t-3\ge 5$.
Notice that this time the residual
$\Res_H(R_{t-3})$ must be empty and so, since
$\mathcal{I}_{\Res_H(R_{t-3})}=\mathcal{O}_{\PP^3}$, we obviously have
$\hh^1({\bf{P}}^3,\mathcal{I}_{\Res_H(R_{t-3})}(t-4)) = 0.$
Hence  thanks to Lemma \ref{castelnuovo} we obtain
$\hh^1({\bf{P}}^3,\mathcal {I}_{Y_t}(t)) = 0$.

We also know that
\begin{equation}\label{stima}
\deg(Y_t)=\deg(X_t)\leq
\binom{t+2}{2}+\binom{t+1}{2}+\binom{t}{2}+\binom{t-1}{2}\leq \binom{t+3}{3}
\end{equation}
where the second inequality is equivalent to
$\binom{t-1}{3}\geq0$, which is
true for any $t\ge 4$. Hence it follows that $Y_t$ has good postulation,
that is ${\A}(t)$ is true.

It is easy to see that  also the scheme $\Res(X_t)$ must be of type
(II) with respect to degree $t-1$. Hence ${\B}(t)$ follows from the
first part of the proof.
\end{proof}

\begin{claim}\label{claim-a}
Fix $t\ge12$.
If $Y_t$ is a zero-dimensional scheme of type (I), then ${\A}(t)$ is
true if ${\B}(t)$ is true.
\end{claim}

\begin{proof}
Since $Y_t$ is of type (I), we can apply Lemma \ref{lemma-a} and we obtain a
specialization $X_t$  such that $\deg(X_t\cap H)=\binom{t+1}{2}$.
Thus, by Lemma \ref{mignon} it follows
$$\hh^0(H,\mathcal {I}_{X_t\cap H}(t))=\hh^1(H,\mathcal {I}_{X_t\cap H}(t))=0.$$

Then, thanks to Lemma \ref{castelnuovo}, it follows, for $i=0,1$,
$$\hh^i({\bf {P}}^3,\mathcal {I}_{X_t}(t))= \hh^i({\bf {P}}^3,\mathcal {I}_{\Res_H(X_t)}(t-1)).$$
Thus in order to prove that the
scheme $X_t$ has good postulation in degree $t$,
it is sufficient to check the good postulation of $\Res_H(X_t)$ in degree
$t-1$.
\end{proof}

\begin{claim}\label{claim-b}
If ${\A}(t-1)$ and ${\C}(t)$ are true, then ${\B}(t)$ is true.
\end{claim}

\begin{proof}
Recall that we write $\Res_H(X_t)=Y_{t-1}\cup Z_{t-1}$,
where $Z_{t-1}$ is a union of simple points supported on $H$.

By Lemma \ref{remark2}, to check that the scheme $\Res_H(X_t)$ has
good postulation in degree $t-1$ (i.e.\ ${\B}(t)$), it is sufficient
to check the good postulation of $Y_{t-1}$ in degree $t-1$ (i.e.\
${\A}(t-1)$) and to prove that ${\C}(t)$ is true.
\end{proof}

\begin{claim}\label{claim-e}
If $Y_t$ is of type (I), then ${\B}(t-1)$ implies ${\C}(t)$.
\end{claim}

\begin{proof}
The statement ${\C}(t)$ is true if 
$\hh^0({\bf {P}}^3,\mathcal {I}_{\Res_H(Y_{t-1})}(t-2)) \le \delta_t.$

Note that since $\deg(X_t\cap H)=\binom{t+2}{2}$, we have
$$ 
\deg(\Res_H(X_t))=\deg(Y_{t-1}\cup Z_{t-1})=\alpha_{t-1}+z_{t-1} =
\alpha_t - \binom{t+2}{2},
$$
and thus it follows
$$\delta_t:=\max\left(0,\binom{t+2}{3} -\alpha_{t-1}-z_{t-1} \right)=
\max\left(0,\binom{t+3}{3} -\alpha_{t}\right).$$

Notice that, by \eqref{caldo} we have $\deg(\Res_H(X_{t})\cap H) \le \binom{t+1}{2}$.
Hence, it follows
$$\deg(\Res_H(Y_{t-1}))=\deg(\Res_H(\Res_H(X_{t})))\ge
\alpha_{t}- \binom{t+2}{2}-\binom{t+1}{2}$$
and then, since $\binom{t+2}{2}+\binom{t+1}{2}=\binom{t+3}{3}-\binom{t+1}{3}$, we get
$$ \deg(\Res_H(Y_{t-1}))\ge
\binom{t+1}{3}-\binom{t+3}{3}+\alpha_t \ge \binom{t+1}{3}-\delta_t.$$

So in order to prove ${\C}(t)$ it is enough to prove that
$\Res_H(Y_{t-1})$ has good postulation in degree $t-2$.
\end{proof}

Now we are in position to prove our main result.
In the following diagram we sketch the steps of the proof:

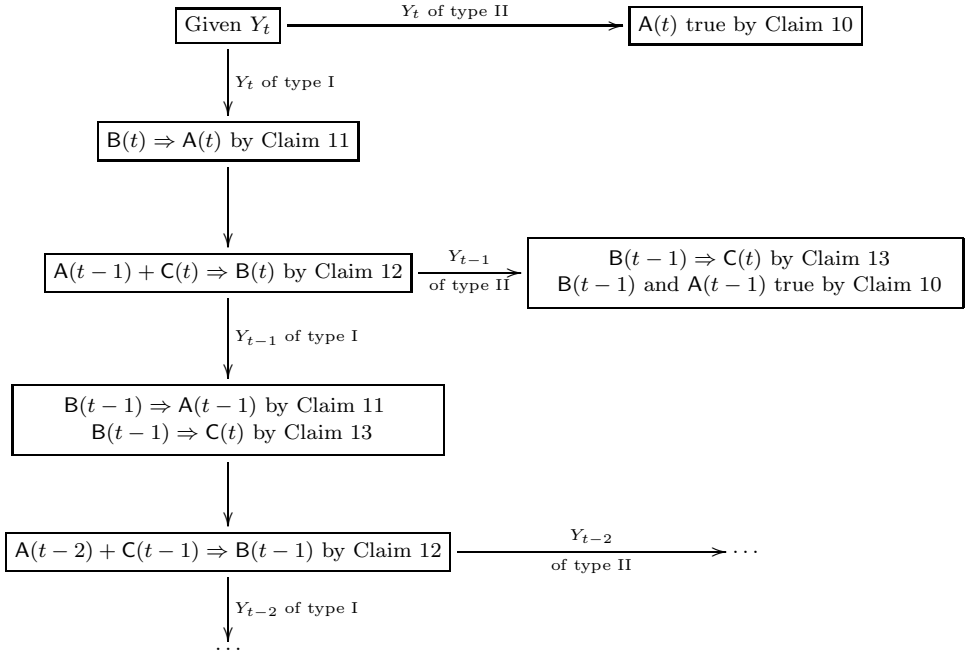
\begin{footnotesize}
\[
    \xymatrix{
      \fbox{Given  $Y_t$} \ar^-{Y_t \textrm{ of type II}}[r] \ar[d]^-{Y_t
        \textrm{ of type I}} &
      \fbox{${\A}(t)$ true by Claim \ref{claim-a0}} \\
      \fbox{${\B}(t)\Rightarrow{\A}(t)$ by Claim \ref{claim-a}} \ar[d] &  \\
      \fbox{${\A}(t-1)+{\C}(t)\Rightarrow{\B}(t)$ by Claim
        \ref{claim-b}} \ar^-{Y_{t-1}}_-{\textrm{of type II}}[r]
      \ar[d]^-{Y_{t-1} \textrm{ of type I}}&
      \fbox{
        \begin{tabular}{c}
          ${\B}(t-1)\Rightarrow{\C}(t)$ by Claim \ref{claim-e}\\
          ${\B}(t-1)$ and ${\A}(t-1)$ true by Claim \ref{claim-a0}
        \end{tabular}}   \\
      \fbox{
        \begin{tabular}{c}
\ \ \           ${\B}(t-1)\Rightarrow{\A}(t-1)$ by Claim \ref{claim-a}\ \ \ \ \ \ \\
          ${\B}(t-1)\Rightarrow{\C}(t)$ by Claim \ref{claim-e}
        \end{tabular}}
      \ar[d] &  \\
      \fbox{${\A}(t-2)+{\C}(t-1)\Rightarrow{\B}(t-1)$ by Claim
        \ref{claim-b}} \ar^-{Y_{t-2}}_-{\textrm{of type II}}[r]
      \ar[d]^-{Y_{t-2} \textrm{ of type I}}&
      {\ldots}\\
      \ldots &
}
\]
\end{footnotesize}

\begin{proof}[Proof of Theorem \ref{i1}]
Fix an integer $d \ge 41$ and a plane $H \subset {\bf {P}}^3$.
For all non negative integers $d,x,y,z,w$, set
$$\epsilon (d,x,y,z):= \binom{d+3}{3} - 20x-10y-4z.$$

Notice that $\epsilon (d,x,y,z+1) =\epsilon (d,x,y,z) -4$,
$\epsilon (d,x,y+1,0)=\epsilon (d,x,y,0)-10$ and
$\epsilon (d,x+1,0,0) = \epsilon (d,x,0,0)-20$.
Hence to prove our statement for all triples $(x,y,z)$ it
is sufficient to check it for all triples $(x,y,z)$ such
that $-19 \le \epsilon (d,x,y,z) \le 3$.

We fix any such triple and
a general union $Y$ of $x$ 4-points, $y$ 3-points and $z$
2-points. We also set $\epsilon=\epsilon (d,x,y,z)$.
We want to prove that $Y$ has good postulation.

Notice that
\begin{equation}\label{numcomponenti}
x+y+z \ge \left\lceil
  \frac{1}{20}\left(\binom{d+3}{3}-3\right)\right\rceil\ge
\frac{1}{20}\binom{d+3}{3}-\frac{3}{20},
\end{equation}
i.e. the scheme $Y$ has at least
$\lceil \frac{1}{20}(\binom{d+3}{3}-3)\rceil$ connected components.

Now we proceed by induction following the steps sketched in the diagram
above. Set $Y_d=Y$.
We can assume by generality that $\deg(Y_d\cap H)\le\binom{d+2}{2}$,
hence we can apply Lemma \ref{lemma-a}, thus specializing the scheme $Y_d$ to a scheme $X_d$.
If $Y_d$ is of type (II), then we conclude by Claim \ref{claim-a0},
since $d\ge 13$.

Hence  we can assume that $Y_d$ is of type (I), and so, since $d\ge
12$, by Claim \ref{claim-a} it is enough to check that the scheme
$\Res_H(X_d)$ has good postulation in degree $d-1$.
Now we write $\Res_H(X_d)=Y_{d-1}\cup Z_{d-1}$,
where $Y_{d-1}$ is the union of all unreduced components of
$\Res_H(X_d)$ and $Z_{d-1}=\Res_H(X_d)\setminus Y_{d-1}$.

By Claim \ref{claim-b}, it is enough to prove that
${\A}(d-1)$ and ${\C}(d)$ are true.

Notice that by \eqref{caldo} we get
$\deg(Y_{d-1}\cap H)\leq \deg(\Res_H(X_d)\cap H)\le \binom{d+1}{2}$.
Hence $Y_{d-1}$ satisfies condition \eqref{condizione} in degree $d-1$, then
we can apply again Lemma \ref{lemma-a}.
We have now two alternatives: either $Y_{d-1}$ is of type (I) or of
type (II).
In both cases, we note that by Claim \ref{claim-e} the statement ${\C}(d)$
follows from ${\B}(d-1)$, since $Y_d$ is of type (I).

Now assume that $Y_{d-1}$ is of type (II). Then by Claim
\ref{claim-a0}, since $d-1\ge13$ we know that ${\B}(d-1)$ and
${\A}(d-1)$ are true and this concludes the proof.
It remains to consider the case $Y_{d-1}$ of type (I).
We apply again Claim \ref{claim-b} and we go on iterating the same
steps.

\smallskip

Now in order to prove our statement we need to show that the
number of steps in the procedure described above is finite.
Moreover we need to check that every time we apply Claim \ref{claim-a0} we are in degree
$\geq13$ and every time we apply Claim \ref{claim-a} we are in degree
$\geq12$.

In order to satisfy all these requirements it is enough to show that in a
finite number of steps we arrive at a scheme of type (II)
in degree $\ge13$.

Recall that we denote, for any integer $t$,
by $X_t$ the specialization described in Lemma \ref{lemma-a},
we write $\Res_H(X_t)=Y_{t-1}\cup Z_{t-1}$ and we set
$z_{t}:= \sharp (Z_{t})$ and
$\alpha_{t}:=\deg(Y_{t})=\deg(X_{t})$.

Now we want to estimate the number of simple points we obtain
iterating the steps above.

Let us assume that starting from the scheme $Y_d$ we arrive in $w$
steps at a scheme $X_{d-w}$ in such a way that the case (II) never
occurs.
Assume also that in these $w$
steps we apply $\gamma$ times Lemma \ref{alehir} with respect to sequences
of type $(1,10,6,3), (1,6,3)$ or $(1,3)$.
Since $g \le 2$, by Lemma \ref{lemma-a}, we have $\gamma\le 2w$. Notice also that the scheme $X_{d-w}$
does not contain simple points, hence it contains at most
$\frac{1}{3}\deg(X_{d-w})=\frac{1}{3}\alpha_{d-w}$ irreducible components.
Hence it follows that
\begin{equation}\label{stima}
\sum _{t=d-w}^{d-1} z_{t}\ge x+y+z-2w-\frac{\alpha_{d-w}}{3}.
\end{equation}
Notice also that
\begin{equation}\label{alfa}
\alpha_{d-w}=\binom{d-w+3}{3}-\epsilon-\sum _{t=d-w}^{d-1} z_{t}.
\end{equation}
and so \eqref{stima} implies
\begin{equation}\label{altra-stima}
\sum _{t=d-w}^{d-1} z_{t}\ge
\frac{3}{2}(x+y+z)-3w-\frac{1}{2}\binom{d-w+3}{3}+\frac{\epsilon}{2}.
\end{equation}

Moreover, setting $w=d$ in \eqref{altra-stima} and
using \eqref{numcomponenti}, we get
\begin{equation}\label{stima-tre}
\sum_{t=0}^{d-1} z_{t} \ge
\frac{3}{40}\binom{d+3}{3}-3d-10-\frac{9}{40}.
\end{equation}

\smallskip

Assume now that $v$ is the maximal integer $w$ such
that for $w$ steps the case (II) is not verified.
Now we prove that we must have $v\le d$.
Indeed the assumption $v>d$ would imply that after $d$ steps we obtain
as a residual a scheme $X_0$ of positive degree (at least $3$).
Hence we have
$$3\le \deg(X_{0})= \alpha_0= -\epsilon+1 -\sum_{t=0}^{d-1} z_{t}\le 20
-\sum_{t=0}^{d-1} z_{t}\le20 $$
 and so $\sum_{t=0}^{d-1} z_{t}\le 17$,
which contradicts \eqref{stima-tre} for $d\geq 17$.

Then let us assume $v\le d$.
Now we want to prove that $d-v\ge13$, and this will conclude the proof.

From the assumption $d \ge 41$ we easily get the following inequality
\begin{equation}\label{bah}
 \frac{1}{20}\binom{d+3}{3}-\frac{3}{20} -19 \ge  2(d-13)
 +\binom{16}{3},
\end{equation}
and then, using \eqref{stima} and \eqref{numcomponenti}, it follows
$$0\le\deg(X_{d-v})= \deg(X_d)-
\left(\binom{d+3}{3}-\binom{d+3-v}{3}\right)-\sum_{t=d-v}^{d-1}
z_{t}\le$$
$$\le -\epsilon+\binom{d+3-v}{3} -
   \frac{1}{20}\binom{d+3}{3} +\frac{3}{20}
 +2v+\frac{\deg(X_{d-v})}{3}. $$
From this we get, by using $\epsilon\ge-19$, $\deg(X_{d-v})\ge0$ and inequality \eqref{bah}
$$0\le -\epsilon+\binom{d+3-v}{3} - \frac{1}{20}\binom{d+3}{3} +\frac{3}{20}
 +2v-\frac{2}{3}\deg(X_{d-v})\le $$
$$\le 19+\binom{d+3-v}{3} - \frac{1}{20}\binom{d+3}{3} +\frac{3}{20}
+2v\le$$
$$\le \binom{d+3-v}{3} +2 (v-d+13)-\binom{16}{3}=:f(d-v).$$

It is easy to  see that $f(d-v)$ is a nondecreasing
function in the interval $d-v\geq 0$, such that $f(13)=0$.
Hence since $f(d-v)\geq 0$, it follows that $d-v\geq 13$, as we wanted.
This concludes the proof of the theorem.
\end{proof}

\section{Cases of low degree: computer aided proofs}

Here we give  the results concerning the cases of low degree, that we
obtained via numerical computations.
\begin{theorem}\label{low}
Assume that $\mathbf{K}$ is an algebraically closed field of characteristic $0$.
Fix non-negative integers $d, x, y, z$ such that $9\le d \le13$. Let $Y
\subset \PP^3$
be a general union of $x$ 4-points, $y$ 3-points and $z$
2-points. Then $Y$ has good postulation.
\end{theorem}

The proof is computer aided and uses the program Macaulay2 \cite{GS}.
Basically we have to check that some matrices, randomly chosen, have
maximal rank. For similar computations see also \cite{partpolint}.

With the same tools, it is not difficult to check that Theorem
\ref{low} is false for $d\le 8$.
For example if we consider  $d=8$ and
$Y$ given by $9$ quartuple points, we expect that there is no
hypersurfaces of degree $8$ through $Y$, but we find that one such
hypersurface exists, since the rank of the corresponding matrix is not
maximal.
Other counterexamples we have found are $8$ quartuple points and $1$
triple point,
$8$ quartuple points and $1$ or $2$ double points, $7$ quartuple
points and $2$ triple points and $1$ double point.
Some of this cases are explained in \cite[Example
7.7]{laface-ugaglia}.

On the other hand, in order to prove Theorem \ref{low} one has to
check a huge number of cases. As an example we list below the
Macaulay2 script which concerns the case $d=12$. Running the
script the computer checks more than 3000 cases, without founding
exceptions. Clearly, with the same method, it is possible to check
the remaining cases $14\le d\le40$. We did not perform this
computation because they need too long time.

Notice that these computations are performed in  characteristic
$31991$, and the result follows in characteristic zero too.
Indeed an integer matrix has maximal rank in
characteristic zero,
if it has maximal rank in positive characteristic.
Furthermore Theorem \ref{low} holds for all positive characteristics with the possible
exception of a finite number of values of the characteristic.

\begin{footnotesize}
\begin{verbatim}
----------------------------------
KK=ZZ/31991;
E=KK[e_0..e_3];
d=12 --degree
N=binomial(d+3,3)

f=ideal(e_0..e_3);
fd=f^d;
T=gens gb(fd)

J=jacobian(T); Jd=J;
--matrix of first derivatives
JJ=jacobian(J);
Jt=submatrix(JJ,{0,1,2,3,5,6,7,10,11,15},{0..N-1});
--matrix of second derivatives: we choose the independent columns
JJJ=jacobian(Jt);
Jq=submatrix(JJJ,{0,2,5,7,9,10,11,14,15,17,19,23,27,30,31,34,35,37,38,39},{0..N-1});
--matrix of third derivatives: we choose the independent columns

mat=random(E^1,E^N)*0
h=1;
for z from 0 to ceiling(N/4) do
for y from 0 to ceiling(N/10) do
for x from 0 to ceiling(N/20) do
(
if ((20*x+10*y+4*z>N-4)and(20*x+10*y+4*z<N+20))
then (print(h,x,y,z), h=h+1,
mat=random(E^1,E^N)*0,
for i from 1 to z do (q=random(E^1,E^4),mat=(mat||sub(Jd,q))),
for i from 1 to y do (q=random(E^1,E^4),mat=(mat||sub(Jt,q))),
for i from 1 to x do (q=random(E^1,E^4),mat=(mat||sub(Jq,q))),
r=rank mat,
if ((20*x+10*y+4*z<N+1)and(r!=20*x+10*y+4*z)) then (print (x,y,z,20*x+10*y+4*z,r)),
if ((20*x+10*y+4*z>N)and(r!=N)) then (print (x,y,z,N,r))))
----------------------------------
\end{verbatim}
\end{footnotesize}

\subsection{The higher multiplicity cases}

It is not difficult to modify the script above in order to perform
some numerical experiments related to the higher multiplicity cases.
Here we list some results obtained 
for schemes of multiplicity $5$. We
denote by $c_i$ the number of $i$-points we consider.

$m=5$, $d=8$:

\begin{center}
\begin{tabular}{|c|c|}
\hline
$(c_5,c_4,c_3,c_2)$ & good postulation\\
\hline
$(5,1,0,0)$& yes \\
$(4,2,0,0)$& no \\
$(3,3,0,0)$& no \\
$(3,4,0,0)$& yes \\
$(2,5,0,0)$& no \\
$(2,6,0,0)$& yes \\
$(1,7,0,0)$& yes \\
$(0,9,0,0)$& no \\
\hline
\end{tabular}
\end{center}

$m=5$, $d=9$:
\begin{center}
\begin{tabular}{|c|c|}
\hline
$(c_5,c_4,c_3,c_2)$ & good postulation\\
\hline
$( 7, 0, 0, 0)$& yes\\
$( 6, 2, 0, 0)$& yes\\
$(5, 3, 0, 0)$& yes\\
$( 4, 5, 0, 0)$& yes\\
$(3,6,0,0)$& no\\
$(3, 7, 0, 0)$& yes\\
$(6,0,1,0)$& no\\
$(6,0,2,0)$& yes\\
\hline
\end{tabular}
\end{center}

$m=5$, $d=10$:
\begin{center}
\begin{tabular}{|c|c|}
\hline
$(c_5,c_4,c_3,c_2)$ & good postulation\\
\hline
$(9, 0, 0, 0)$& no\\
$(8, 1, 0, 0)$& no\\
$(7, 2, 0, 0)$& no\\
$(8, 2, 0, 0)$& yes\\
$(7, 3, 0, 0)$& yes\\
$(6, 4, 0, 0)$& yes\\
$(6, 5, 0, 0)$& yes\\
$(8, 0, 1, 0)$& no\\
\hline
\end{tabular}
\end{center}

\medskip

\providecommand{\bysame}{\leavevmode\hbox to3em{\hrulefill}\thinspace}


\begin{thebibliography}{99}

\bibitem{AHinv} James Alexander and Andr\'e Hirschowitz,
La m\'ethode d'Horace \'eclat\'ee: application \`a l'interpolation en
degr\'ee quatre, Invent. Math. 107 (1992), 585--602.

\bibitem{AHjag} James Alexander and Andr\'e Hirschowitz, Polynomial
  interpolation in several variables, J. Alg. Geom. 4 (1995), no.2, 201--222.

\bibitem{AH}
James Alexander and Andr\'e Hirschowitz, An asymptotic vanishing
theorem for generic unions of multiple points, Invent. Math.
{140} (2000), no. 2, 303--325.

\bibitem{ballico-tripli} Edoardo Ballico. On the postulation of general union of double points and triple points
in $\bf{P}^3$, preprint 2007.

\bibitem{BO} Maria Chiara Brambilla and Giorgio Ottaviani, On the
  Alexander-Hirschowitz Theorem, J. Pure Appl. Algebra 212
  (2008), no. 5, 1229-1251.

\bibitem{partpolint} Maria Chiara Brambilla and Giorgio Ottaviani, On
  partial polynomial interpolation, preprint 2007,
 arXiv:math/0705.4448


\bibitem{ciliberto} Ciro Ciliberto, Geometrical aspects of polynomial
  interpolation in more variables and
of Waring's problem, ECM Vol. I (Barcelona, 2000),
 Progr. Math., 201, Birkh\"auser, Basel, 2001, 289-316.

\bibitem{dumnicki} Marcin Dumnicki, Regularity and non-emptyness of
  linear systems in $\PP^n$, preprint 2008, arXiv:0802.0925.

\bibitem{GS} Daniel Grayson and Michael Stillman,
          Macaulay 2, a software system for research
                   in algebraic geometry, available at http://www.math.uiuc.edu/Macaulay2/

\bibitem{m} Thierry Mignon, Syst\`{e}mes de courbes planes \`{a} singularit\'{e}s impos\'{e}es: le cas de multiplicit\'{e}s
inf\'{e}rieurs ou \'{e}gales \`{a} quatre, J. Pure Appl. Algebra 151 (2000), no. 2, 173--195.


\bibitem{laface-ugaglia} Antonio Laface and Luca Ugaglia,
{On a class of special linear systems of $\PP^3$},
{Trans. Amer. Math. Soc.}, {358} (2006), no. 12, 5485--5500 (electronic).

\bibitem{volder-laface} Cindy De Volder and Antonio Laface,
{On linear systems of $\PP^3$ through multiple points},
{J. Algebra}, 310 (2007), no. 1, 207-217.

\end{thebibliography}
\end{document}